\documentclass{birkmult} \usepackage{jms-bm}
\begin{document}



%


\newcommand{\intersect}{\cap}
\newcommand{\union}{\cup}
\newcommand{\R}{\mathbb{R}}
\newcommand{\E}{\mathbb{E}}
\let\Section=\S
\renewcommand{\S}{\mathbb{S}}
\renewcommand{\H}{\mathbb{H}}
\renewcommand{\O}{\mathcal{O}}
\newcommand{\vH}{\mathbf{H}} 
\newcommand{\vx}{\mathbf{x}} 
\newcommand{\vA}{\mathbf{A}} 
\newcommand{\area}{\operatorname{Area}}
\newcommand{\vol}{\operatorname{Vol}}
\newcommand{\len}{\operatorname{Len}}
\newcommand{\Link}{\operatorname{Link}}
\newcommand{\Star}{\operatorname{Star}}
\renewcommand{\id}{\operatorname{Id}}
\renewcommand{\d}{\partial}
\newcommand{\del}{\nabla}
\newcommand{\eps}{\varepsilon}
\newcommand{\<}{\langle}
\renewcommand{\>}{\rangle}
\newcommand{\cross}{\times}
\newcommand{\half}{{\frac12}}
\newcommand{\isom}{\cong}
\newcommand{\ident}{\equiv}

\def\Gauss/{Gauss}



\title{Curvatures of Smooth and Discrete Surfaces}

\author{John M. Sullivan}
\email{sullivan@math.tu-berlin.de}
\address{\tuaddr32}


\def\incgrh#1#2{\includegraphics[height=#2]{figs/#1}}
\def\figr#1{Figure~\ref{fig:#1}}

\begin{abstract}
We discuss notions of Gauss curvature and mean curvature
for polyhedral surfaces.  The discretizations are guided
by the principle of preserving integral relations for
curvatures, like the Gauss/Bonnet theorem and
the mean-curvature force balance equation.
\end{abstract}

\maketitle

The curvatures of a smooth surface are local measures
of its shape.  Here we consider analogous quantities for discrete
surfaces, meaning triangulated polyhedral surfaces.
Often the most useful analogs are those which preserve
integral relations for curvature, like the \Gauss//Bonnet theorem
or the force balance equation for mean curvature.
For simplicity, we usually restrict our attention to surfaces
in euclidean three-space~$\E^3$, although some of the results
generalize to other ambient manifolds of arbitrary dimension.

This article is intended as background for some of the related
contributions to this volume.
Much of the material here is not new; some is even quite old.
Although some references are given, no attempt has been made to
give a comprehensive bibliography or a full picture of the
historical development of the ideas.

\section{Smooth curves, framings and integral curvature relations}

A companion article~\cite{Sull-FTC} in this volume investigates curves
of finite total curvature.  This class includes both
smooth and polygonal curves, and allows a unified treatment
of curvature.  Here we briefly review the theory of smooth curves
from the point of view we will later adopt for surfaces.

The \index{curvature (of curve)}{curvatures of a smooth curve}~$\gamma$ (which we usually assume is
parametrized by its arclength~$s$) are the local properties of its shape,
invariant under euclidean motions.
The only first-order information is the tangent line;
since all lines in space are equivalent, there are no first-order invariants.
Second-order information (again, independent of parametrization)
is given by the osculating circle; the one corresponding invariant
is its curvature $\kappa=1/r$.

(For a plane curve given as a graph $y=f(x)$ let us contrast the notions
of curvature~$\kappa$ and second derivative~$f''$.  At a point~$p$ on the
curve, we can find either one by translating~$p$ to the origin, transforming
so the curve is horizontal there, and then looking at the second-order
behavior.  The difference is that for curvature, the
transformation is a euclidean rotation, while for second derivative,
it is a shear $(x,y)\mapsto (x,y-ax)$.  A~parabola has constant
second derivative $f''$ because it looks the same at any two points
after a shear.  A circle, on the other hand, has constant curvature
because it looks the same at any two points after a rotation.)

A plane curve is completely determined (up to rigid motion)
by its (signed) curvature~$\kappa(s)$ as a function of arclength $s$.
For a space curve, however, we need to look at the third-order
invariants; these are the torsion~$\tau$ and the derivative~$\kappa'$,
but the latter of course gives no new information.
Curvature and torsion now form a complete set of invariants:
a space curve is determined by~$\kappa(s)$ and~$\tau(s)$.

Generically speaking, while second-order curvatures usually suffice to
determine a hypersurface (of codimension~$1$),
higher-order invariants are needed for higher codimension.
For curves in~$\E^d$, for instance, we need $d-1$ generalized
curvatures, of order up to~$d$, to characterize the shape.

Let us examine the case of space curves $\gamma\subset\E^3$ in more detail.
At every point $p\in\gamma$ we have a splitting of the tangent space~$T_p\E^3$
into the tangent line~$T_p\gamma$ and the normal plane.
A framing along~$\gamma$ is a smooth choice
of a unit normal vector~$N_1$, which is then completed to the oriented
orthonormal frame $(T,N_1,N_2)$ for~$\E^3$, where $N_2=T\cross N_1$. 
Taking the derivative with respect to arclength, we get a skew-symmetric
matrix (an infinitesimal rotation) that describes how the frame changes:
$${\begin{pmatrix} T \\ N_1 \\ N_2 \end{pmatrix}\!}'\,
= \begin{pmatrix} 0 & \kappa_1 & \kappa_2 \\
               -\kappa_1 & 0 & \tau \\
               -\kappa_2 & -\tau & 0 \end{pmatrix}
\begin{pmatrix} T \\ N_1 \\ N_2 \end{pmatrix}\!.
$$
Here, $T'(s)=\sum\kappa_i N_i$ is the \emph{curvature vector} of~$\gamma$,
while~$\tau$ measures the twisting of the chosen orthonormal frame.

If we fix $N_1$ at some basepoint along~$\gamma$, then one natural framing
is the \dfn{parallel frame} or \emph{Bishop frame}~\cite{Bishop-frame}%
\index{Bishop frame|seeonly{parallel frame}}
defined by the condition~$\tau=0$.  Equivalently, the vectors~$N_i$
are parallel-transported along~$\gamma$ from the basepoint, using the
Riemannian connection on the normal bundle induced by the immersion in~$\E^3$.
One should note that this is not necessarily a closed framing along
a closed loop~$\gamma$; when we return to the basepoint, the vector~$N_1$ has
been rotated through an angle called the \dfn{writhe} of~$\gamma$.

Other framings are also often useful.
For instance, if~$\gamma$ lies on a surface~$M$ with unit normal~$\nu$,
it is natural to choose $N_1=\nu$.  Then $N_2=\eta := T\cross\nu$
is called the \emph{cornormal} vector, and $(T,\nu,\eta)$ is the
\dfn{Darboux frame} (adapted to $\gamma\subset M\subset\E^3$).
The curvature vector of~$\gamma$ decomposes into parts tangent and normal
to~$M$ as $T'= \kappa_n \nu + \kappa_g \eta$.
Here, $\kappa_n=-\nu'\cdot T$ measures
the \emph{normal curvature} of~$M$ in the direction~$T$, and is independent of~$\gamma$, 
while~$\kappa_g$, the \emph{geodesic curvature} of~$\gamma$ in~$M$,
is an intrinsic notion, showing how~$\gamma$ sits in~$M$,
and is unchanged if we isometrically deform the immersion of~$M$ into space.

When the curvature vector of~$\gamma$ never vanishes, 
we can write it as $T'=\kappa N$, where~$N$~is
a unit vector, the \emph{principal normal}, and $\kappa>0$.
This yields the orthonormal \dfn{Frenet frame} $(T,N,B)$,
whose twisting~$\tau$ is the \emph{torsion} of~$\gamma$.

The \emph{total curvature} of a smooth curve is $\int\!\kappa\,ds$.
In~\cite{Sull-FTC} we review a number of standard results:
For closed curves, the total curvature is at least~$2\pi$
(Fenchel) and for knotted space curves the total curvature
is at least $4\pi$ (F\'ary/Milnor).  For plane curves,
we can consider instead the signed curvature, and find
that $\int\!\kappa\,ds$ is always an integral multiple of~$2\pi$.
Suppose (following Milnor) we define the total curvature of a
polygonal curve simply to be the sum of the turning angles
at the vertices.  Then, as we explain in~\cite{Sull-FTC},
all these theorems on total curvature
remain true.  Our goal, when defining curvatures for polyhedral
surfaces, will be to ensure that similar integral relations
remain true.

\section{Curvatures of smooth surfaces}

Given a (two-dimensional, oriented) surface~$M$ smoothly immersed in~$\E^3$,
we understand its local shape by looking at the \ix{Gauss map}
$\nu: M\to\S^2$ given by the unit normal vector $\nu=\nu_p$
at each point $p\in M$.
The derivative of the \Gauss/ map at~$p$ is a linear map
from~$T_pM$ to~$T_{\nu_p}\S^2$.
Since these spaces are naturally identified, being parallel planes in~$\E^3$,
we can view the derivative as an endomorphism $-S_p: T_pM \to T_pM$.
The map~$S_p$ is called the \ix{shape operator} (or Weingarten map).

The shape operator is the complete second-order invariant (or curvature)
which determines the original surface~$M$.  (This statement has been left
intentionally a bit vague, since without a standard parametrization
like arclength, it is not quite clear how one should specify such an
operator along an unknown surface.)  Usually, however, it is more convenient
to work not with the operator~$S_p$ but instead with scalar quantities.
Its eigenvalues~$\kappa_1$ and~$\kappa_2$ are called the
\emph{principal curvatures}, and (since they cannot be globally distinguished)
it is their symmetric functions which have the most geometric meaning.

We define the \dfn{Gauss curvature} $K:=\kappa_1\kappa_2$
as the determinant of $S_p$ and the \emph{mean curvature}
$H:=\kappa_1+\kappa_2$ as its trace.  Note that the sign of
$H$ depends on the choice of unit normal $\nu$, so often
it is more natural to work with the
\emph{vector mean curvature}\index{mean curvature!vector}
(or mean curvature vector) $\vH:=H\nu$.  Furthermore,
some authors use the opposite sign on~$S_p$ and thus~$\vH$,
and many use $H=(\kappa_1+\kappa_2)/2$, justifying the
name ``mean'' curvature.  Our conventions mean that the mean
curvature vector for a convex surface points inwards (like the curvature
vector for a circle).  For a unit sphere oriented with inward normal,
the \Gauss/ map~$\nu$ is the antipodal map, $S_p=I$, and $H\ident2$.

The \Gauss/ curvature is an intrinsic notion, depending only on
the pullback metric on the surface~$M$, and not on the immersion into space.
That is, $K$ is unchanged by bending the surface without stretching it.
For instance, a developable surface like a cylinder or cone has $K\ident0$
because it is obtained by bending a flat plane.  One intrinsic characterization
of $K(p)$ is obtained by comparing the circumferences $C_\eps$ of (intrinsic)
$\eps$-balls around~$p$ to the value $2\pi\eps$ in~$\E^2$.  We get
$$\frac{C_\eps}{2\pi\eps} = 1 - \frac{\eps^2}6 K(p) + \O(\eps^3).$$ 

{Mean curvature}, on the other hand,
is certainly not intrinsic, but it has a nice variational
interpretation.  Consider a variation vectorfield~$\xi$ on~$M$;
for simplicity assume~$\xi$ is compactly supported away from any boundary.
Then $H=-\delta\area/\delta\vol$ in the sense that
$$\delta_\xi\vol = \int \xi\cdot\nu\,dA, \qquad
  \delta_\xi\area = -\int \xi\cdot H\nu\,dA.$$
With respect to the $L^2$ inner product
$\<\xi,\eta\> := \int \xi_p\cdot \eta_p\,dA$
on vectorfields, the vector mean curvature
is thus the negative gradient of the area functional, often called the
first variation of area: $\vH=-\del\area$.
(Similarly, the negative gradient of length for a curve is
its curvature vector~$\kappa N$.)

Just as $\kappa$ is the geometric version of second derivative
for curves, mean curvature is the geometric version of the Laplacian $\Delta$.
Indeed, if a surface $M$ is written locally near a point~$p$
as the graph of a height function~$f$ over its tangent plane $T_p M$,
then $H(p)=\Delta f$.
Alternatively, we can write $\vH = \del_M\cdot\nu = \Delta_M \vx$, where
$\vx$ is the position vector in $\E^3$ and $\Delta_M$ is the
\index{Laplace--Beltrami operator}Laplace/Beltrami operator,
the intrinsic surface Laplacian.

We can flow a curve or surface to reduce its length or area,
by following the gradient vectorfield $\kappa N$ or $H\nu$; the
resulting parabolic (heat) flow is slightly nonlinear
in a natural geometric way.  This so-called \emph{mean-curvature flow}
has been extensively studied as a geometric smoothing flow. 
(See, among many others, \cite{GagHam-csf,Gra-csf} for the
curve-shortening flow
and \cite{Bra-mcf,Huisken-mcf,Ilmanen-mcf,Ecker-mcf,White-mcf}
for higher dimensions.)

\section{Integral curvature relations for surfaces}

For surfaces, we will consider various integral curvature relations
that relate area integrals over a region $D\subset M$
to arclength integrals over the boundary $\gamma:=\d D$.
First, the \Gauss//Bonnet theorem\index{Gauss Bonnet@Gauss/Bonnet theorem} says that, when $D$ is a disk,
$$2\pi-\iint_D K\,dA = \oint_\gamma \kappa_g\,ds
   = \oint_\gamma T'\cdot\eta\,ds = -\oint_\gamma\eta'\cdot d\vx.$$
Here, $d\vx = T\,ds$ is the vector line element along~$\gamma$,
and $\eta=T\times\nu$ is again the cornormal.
In particular, this theorem implies that the total \Gauss/ curvature
of~$D$ depends only on a collar neighborhood of $\gamma$:
if we make any modification to~$D$ supported away from the boundary,
the total curvature is unchanged (as long as~$D$ remains topologically a disk).
We will extend the notion of (total) \Gauss/ curvature from smooth surfaces
to more general surfaces (in particular polyhedral surfaces) by
requiring that this property remain true.

Our other integral relations are all proved by Stokes's Theorem, and thus
require only that~$\gamma$ be the boundary of~$D$ in a homological sense;
for these $D$ need not be a disk.  First consider the
\dfn{vector area}
$$\vA_\gamma := \half\oint_\gamma \vx\cross d\vx
= \half\oint_\gamma \vx\cross T\,ds = \iint_D \nu\,dA.$$
The right-hand side represents the total vector area of any surface
spanning $\gamma$, and the relation shows this to depend only on~$\gamma$
(and this time not even on a collar neighborhood).  The integrand on
the left-hand side depends on a choice of origin for the coordinates,
but because we integrate over a closed loop, the integral is independent
of this choice.  Both sides of this vector area formula can be interpreted
directly for a polyhedral surface, and the equation remains true in that case.
We note also that this vector area~$\vA_\gamma$ is one of the quantities
preserved when~$\gamma$ evolves under the Hasimoto\index{Hashimoto flow/surface}
or smoke-ring flow
$\dot\gamma = \kappa B$.  (Compare \cite{LanPer,PSW,HoffmannSRF}.)

A simple application of the fundamental theorem of calculus
to the tangent vector of a curve~$\gamma$ from~$p$ to~$q$ shows that
$$T(q)-T(p) = \int_p^q T'(s)\,ds = \int_p^q \kappa N\,ds.$$
This can be viewed as a balance between elastic tension forces trying to
shrink the curve, and sideways forces holding it in place.  It is the
key step in verifying that the vector curvature $\kappa N$ is the
first variation of length.

The analog for a surface patch~$D$ is the
mean-curvature \ix{force balance} equation
$$\oint_\gamma \eta\,ds = -\oint_\gamma \nu\cross d\vx
 = \iint_D H\nu\,dA = \iint_D \vH\,dA.$$
Again this represents a balance between surface tension forces acting
in the conormal direction along the boundary of $D$ and what can be
considered as pressure forces (espcially in the case of constant~$H$)
acting normally across~$D$.  We will use this equation to develop our
analog of mean curvature for discrete surfaces.

The force balance equation can be seen to arise from the translational
invariance of curvatures.  It has been important for studying surfaces
of constant mean curvature; see for instance \cite{KKS, Kus-bub, triund}.
The rotational analog is the following \ix{torque balance}:
$$\oint_\gamma \vx\cross\eta\,ds = \oint_\gamma\vx\cross(\nu\cross d\vx)
 = \iint_D H(\vx\cross\nu)\,dA = \iint_D \vx\cross\vH\,dA.$$
Somewhat related is the following equation:
$$\oint_\gamma \vx\cdot\eta\,ds = \oint_\gamma\vx\cdot (\nu\cross d\vx)
 = \iint_D(\vH\cdot\vx - 2)\,dA.$$
It gives, for example, an interesting
expression for the area of a minimal ($H\ident0$) surface.

\section{Discrete surfaces}
For us, a discrete or polyhedral
surface\index{surface!polyhedral}
$M\subset\E^3$ will mean
a triangulated surface with a continuous map into space, linear on each
triangle.  In more detail, we start with an abstract combinatorial
\ix{triangulation}---a simplicial complex---representing
a $2$-manifold with boundary.
We then pick positions $p\in\E^3$ for all the vertices,
which uniquely determine a linear map on each triangle;
these maps fit together to form the polyhedral surface.

The union of all triangles containing a vertex~$p$
is called $\Star(p)$, the \emph{star} of~$p$.
Similarly, the union of the two triangles
containing an edge~$e$ is~$\Star(e)$.

\subsection{\Gauss/ curvature}
\index{Gauss curvature!discrete}
It is well known how the notion of \Gauss/ curvature extends to such
discrete surfaces~$M$.
(Banchoff~\cite{Banchoff-jdg,Banchoff-monthly} was probably
the first to discuss this in detail, though he notes that Hilbert and
Cohn-Vossen~\cite[\Section 29]{HCV-AG}
had already used a polyhedral analog to motivate
the intrinsic nature of \Gauss/ curvature.)
Any two adjacent triangles (or, more generally,
any simply connected region in~$M$ not including any vertices) can
be flattened---developed isometrically into the plane.  Thus the \Gauss/
curvature is supported on the vertices $p\in M$.  In fact, to
keep the \Gauss//Bonnet theorem true, we must take
$$\iint_D K\,dA := \sum_{p\in D} K_p, \qquad
\text{with} \quad K_p := 2\pi - \sum_i\theta_i.$$
Here, the angles~$\theta_i$ are the interior angles at~$p$
of the triangles meeting there, and $K_p$ is often known as the angle
defect at~$p$.  If~$D$ is any neighborhood of~$p$ contained in $\Star(p)$,
then $\oint_{\d D} \kappa_g\,ds=\sum \theta_i$; when the triangles are
acute, this is most easily seen by letting~$\d D$ be the path connecting
their circumcenters and crossing each edge perpendicularly.
Analogous to our intrinsic characterization of \Gauss/ curvature
in the smooth case, note that the circumference of a small $\eps$ ball
around~$p$ here is exactly $2\pi\eps-\eps K_p$.

This version of discrete \Gauss/ curvature
is quite natural, and seems to be the correct analog
when \Gauss/ curvature is used intrinsically.
But the \Gauss/ curvature of a smooth surface in~$\E^3$
also has extrinsic meaning; for instance the total
absolute \Gauss/ curvature is proportional to the average
number of critical points of different height functions.
For such considerations, Brehm
and K{\"u}hnel~\cite{BreKuh} suggest the following: when a vertex~$p$
is extreme on the convex hull of its star, but the star itself
is not a convex cone, then we should think of~$p$ as having both
positive and negative curvature.
We let $K^+_p$ be the curvature at~$p$ of the convex hull,
and set $K^-_p:= K^+_p-K_p\ge0$. Then the absolute curvature
at~$p$ is $K^+_p + K^-_p$, which is greater than $|K_p|$.
This discretization, is of course, also based on preserving
an integral curvature relation---a different one.
(See also~\cite{BanKuh,vDAlb}.)

We can use the same principle as before---preserving
the \Gauss//Bonnet theorem---to define \Gauss/ curvature, as a measure, for 
much more general surfaces.  For instance, on a piecewise
smooth surface, we have ordinary~$K$ within each face, a point
mass (again the angle defect) at each vertex,
and a linear density along each edge, equal to the
difference in the geodesic curvatures of that edge within its
two incident faces.  Indeed, clothes are often designed from
pieces of (intrinsically flat) cloth, joined so that each vertex
is intrinsically flat and thus all the curvature is along the edges;
corners would be unsightly in clothes.

Returning to polyhedral surfaces, we
note that~$K_p$ is clearly an intrinsic notion (as it should be)
depending only on the angles of each triangle and not on the precise
embedding into~$\E^3$.  Sometimes it is useful to have a notion of
combinatorial curvature, independent of \emph{all} geometric information.
Given just a combinatorial triangulation, we can pretend that each
triangle is equilateral with angles $\theta=60^\circ$.
(Such a euclidean metric with cone points at certain vertices exists on the
abstract surface, independent of whether or not it could be embedded in space.
See the survey~\cite{Troyanov-survey}.)

The curvature of this metric, $K_p = \frac\pi3 (6-\deg p)$, can be called the
\emph{combinatorial (\Gauss/) curvature}\index{Gauss curvature!combinatorial}
of the triangulation.  (See~\cite{Thurston-triang,57torus}
for combinatorial applications of this notion.)
In this context, the global form $\sum K_p = 2\pi\chi(M)$
of \Gauss//Bonnet amounts to nothing more than
Euler's formula $\chi=V-E+F$.\index{Euler characteristic}
(We note that Forman has proposed a combinatorial
Ricci curvature~\cite{Forman-combRic}; although for smooth surfaces,
Ricci curvature is \Gauss/ curvature, for discrete surfaces
Forman's combinatorial curvature does not agree with ours,
so he fails to recover the \Gauss//Bonnet theorem.)

Our discrete \Gauss/ curvature~$K_p$ is of course an integrated
quantity.  Sometimes it is desirable to have instead a curvature
density, dividing~$K_p$ by the surface area associated to the
vertex~$p$.  One natural choice is $A_p:= \tfrac13 \area(\Star(p))$,
but this does not always behave nicely for irregular triangulations.
One problem is that, while $K_p$ is intrinsic, depending only on
the cone metric of the surface, $A_p$ depends also on the choice
of which pairs of cone points are connected by triangle edges.
One fully intrinsic notion of the area associated to~$p$ would be
the area of its intrinsic Voronoi cell in the sense of
Bobenko and Springborn~\cite{BobSpr-discrLB};
perhaps this would be the best choice for computing \Gauss/ curvature density.

\subsection{Vector area}
The \ix{vector area} formula
$$\vA_\gamma := \half\oint_\gamma \vx\cross d\vx = \iint_D \nu\,dA$$
needs no special interpretation for discrete surfaces: both sides
of the equation make sense directly, since the surface normal $\nu$ is
well-defined almost everywhere.  However, it is worth interpreting
this formula for the case when $D$ is the star of a vertex $p$. 
More generally, suppose $\gamma$ is any closed curve (smooth or polygonal),
and $D$ is the cone from $p$ to $\gamma$ (the union of all line segments
$pq$ for $q\in\gamma$).  Fixing $\gamma$ and letting $p$ vary, we find that
the volume enclosed by this cone is an affine linear function of $p\in\E^3$,
and thus
$$\vA_p :=  \del_p\vol D = \frac{\vA_{\gamma}}3
         = \frac16\oint_\gamma \vx\cross d\vx$$
is independent of the position of~$p$.
We also note that any such cone $D$ is intrinsically flat except at the
cone point $p$, and that $2\pi-K_p$ is the cone angle at $p$.

\subsection{Mean curvature}
The mean curvature\index{mean curvature!discrete}
of a discrete surface~$M$ is supported along the edges.
If $e$ is an edge, and $e\subset D\subset\Star(e)=T_1\cup T_2$, then we set
$$\vH_e := \iint_D\vH\,dA = \oint_{\d D}\eta\,ds
         = e\cross\nu_1 - e\cross\nu_2 = J_1e - J_2e.$$ 
Here $\nu_i$ is the normal vector to the triangle $T_i$,
and the operator $J_i$ rotates by $90^\circ$ in the plane of that triangle.
Note that $|\vH_e| = 2 \sin(\theta_e/2)\,|e|$, where $\theta_e$
is the exterior dihedral angle along the edge, defined by
$\cos\theta_e = \nu_1\cdot\nu_2$.

No nonplanar discrete surface has $\vH_e=0$ along every edge.
But this discrete mean curvature can cancel out around vertices.  We set
$$2\vH_p := \sum_{e\ni p} \vH_e = \iint_{\Star(p)}\!\!\!\vH\,dA
          = \oint_{\d\Star(p)}\!\!\!\eta\,ds.$$
The area of the discrete surface is a function of the vertex positions;
if we vary only one vertex $p$, we find that
$\del_p\area(M) = -\vH_p$.  This mirrors the variational characterization
of mean curvature for smooth surfaces, and we see that a natural
notion of discrete minimal surfaces\index{minimal surface (discrete)}
is to require $\vH_p\ident0$ for all vertices~\cite{Bra-evol, PinPol}.

Suppose that the vertices adjacent to $p$, in cyclic order,
are $p_1$, \ldots, $p_n$.
Then we can express~$\vA_p$ and~$\vH_p$ explicitly in terms
of these neighbors.  We get
\begin{align*}
3\vA_p &= 3\del_p\!\vol = \iint_{\Star(p)}\!\!\!\nu\,dA
        = \half\oint_{\d\Star(p)}\!\!\!\vx\cross d\vx
             = \half\sum_i p_i\cross p_{i+1}
\end{align*}
and similarly
\begin{align*}
2\vH_p &= \sum\vH_{pp_i} = -2\del_p\area = \sum J_i(p_{i+1}-p_i) \\
       &= \sum_i (\cot\alpha_i + \cot\beta_i)(p-p_i),
\end{align*}
where $\alpha_i$ and $\beta_i$ are the angles opposite edge $pp_i$ in the
two incident triangles.  This latter equation is the
famous ``\ix{cotangent formula}''~\cite{PinPol, Wardetzky}
which also arises naturally in a finite-element discretization of
the Laplacian.

Suppose we change the combinatorics of a discrete surface~$M$ by
introducing a new vertex~$p$ along an existing edge~$e$ and subdividing
the two incident triangles.  Then $\vH_p$ in the new surface equals
the original $\vH_e$, independent of where along~$e$ we place~$p$.
This allows a variational interpretation of~$\vH_e$.

\subsection{Minkowski mixed volumes}
A somewhat different interpretation of mean curvature for convex polyhedra
is found in the context of Minkowski's theory of
mixed volumes.\index{intrinsic/mixed volumes}
(In this simple form, it dates back well before Minkowski,
to Steiner~\cite{Steiner-mv}.)
If $X$ is a smooth convex
body in $\E^3$ and $B_t(X)$ denotes its $t$-neighborhood, then 
its Steiner polynomial\index{Steiner's formula} is:
$$\vol(B_t(X)) = \vol X + t \area \d X + \frac{t^2}2 \int_{\d X} \!\!H\,dA 
                        + \frac{t^3}3 \int_{\d X} \!\!K\,dA.$$
Here, by \Gauss//Bonnet, the last integral is always $4\pi$.

When $X$ is instead a convex polyhedron, we already understand how to
interpret each term except $\int_{\d X} H\,dA$.
The correct replacement for this term, as Steiner discovered,
is $\sum_e \theta_e\,|e|$.
This suggests $H_e := \theta_e\,|e|$ as a notion of total mean curvature
for the edge~$e$.

Note the difference between this formula and our earlier
$|\vH_e| = 2\sin(\theta_e/2)\,|e|$.  Either one can be derived by
replacing the edge $e$ with a sector of a cylinder of length~$|e|$
and arbitrary (small) radius $r$.  We have $\iint\vH\,dA=\vH_e$,
so that
$$\biggl|\iint\vH\,dA\biggr| = |\vH_e|=2\sin(\theta_e/2)\,|e|
   \;\, < \;\, \theta_e\,|e| = H_e =\iint H\,dA. $$
The difference is explained by the fact that one formula integrates
the scalar mean curvature while the other integrates the vector mean curvature.
Again, these two discretizations both arise through preservation of
(different) integral relations for mean curvature.

See~\cite{Sull-FTC} for a more extensive discussion of the analogous
situation for curves: although as we have mentioned,
the sum of the turning angles~$\psi_i$
is often the best notion of total curvature for a polygon,
in certain situations the ``right'' discretization is instead the sum of
$2\sin\psi_i/2$ or $2\tan\psi_i/2$.

The interpretation of curvatures in terms of the mixed volumes or
Steiner polynomial actually works for arbitrary convex surfaces.
(Compare~\cite{Schroeder} in this volume.)
Using geometric measure theory---and a generalized normal
bundle called the normal cycle---one can extend both
\Gauss/ and mean curvature in a similar way to quite general surfaces.
See~\cite{Fed-curvmeas,Fu-curvsuban,CSM-norcyc,CSM-sfm}.

\subsection{Constant mean curvature and Willmore surfaces}
A smooth surface which minimizes area under a volume constraint
has constant mean curvature; the constant $H$ can be understood
as the Lagrange multiplier for the constrained minimization problem.
A discrete surface which minimizes area among surfaces
of fixed combinatorial type and fixed volume will have constant
discrete mean curvature $H$ in the sense that at every vertex,
$\vH_p = H\vA_p$, or equivalently $\del_p\area = -H \del_p\!\vol$. 

In general, of course, the vectors $\vH_p$ and $\vA_p$ are not
even parallel: they give two competing notions of a normal vector to the
discrete surface at the vertex~$p$.
Still,
$$h_p := \frac{\big|\del_p\area\big|}{\big|\del_p\vol\big|} = \frac{\big|\vH_p\big|}{\big|\vA_p\big|} =
\frac{\big|\iint_{\Star(p)}\vH\,dA\big|}{\big|\iint_{\Star(p)}\nu\,dA\big|}$$
gives a better notion of mean curvature density near~$p$ than, say,
the smaller quantity
$$\frac{|\vH_p|}{A_p}
  = \frac{\big|\iint\vH\,dA\big|}{\iint 1\,dA},$$
where $A_p:= \tfrac13 \area(\Star(p))$.

Suppose we want to discretize the elastic bending energy for surfaces,
$\iint H^2\,dA$,
known as the \ix{Willmore energy}.
The discussion above
shows why $\sum_p h_p^2 A_p$ (which was used in~\cite{minimax})
is a better discretization than $\sum_p |\vH_p|^2/A_p$
(used fifteen years ago in~\cite{HKS}).
Several related discretizations are by now built into Brakke's Evolver;
see the discusion in~\cite{Bra-manual}.
Recently, Bobenko~\cite{Bob-willm,Bobenko-circ} has described
a completely different approach to discretizing the Willmore energy,
which respects the M{\"o}bius invariance of the smooth energy.

\subsection{Relation to discrete harmonic maps}
As mentioned above, we can define a \emph{discrete minimal
surface} to be a polyhedral
surface with $\vH_p\ident0$.
An early impetus to the field of discrete differential
geometry was the realization (starting with~\cite{PinPol})
that discrete minimal surfaces are not only critical points
for area (fixing the combinatorics), but also have other
properties similar to those of smooth minimal surfaces.

For instance, in a conformal parameterization of a
smooth minimal surface, the coordinate functions are harmonic.
To interpret this for discrete surfaces, we are led to the
question of when a discrete map should be considered conformal.
In general this is still open.
(Interesting suggestions come from the theory of circle packings,
and this is an area of active research.
See for instance~\cite{Steph-CPbook,BowersHurdal,Bobenko-circ,
KSS,Springborn-delcirc}.)

However, we should certainly agree that the identity map is conformal.
A polyhedral surface $M$ comes with an embedding $\id_M:M\to\E^3$
which we consider as the identity map.  Indeed, we then find
(following~\cite{PinPol}) that $M$ is discrete minimal
if and only if $\id_M$ is \index{harmonic function (discrete)}{discrete harmonic}.
Here a polyhedral map $f:M\to \E^3$ is called \emph{discrete
harmonic} if it is a critical point for the Dirichlet energy,
written as the following sum over the triangles~$T$ of~$M$:
$$E(f) := \sum_T |\del f_T|^2 \area_M(T).$$
We can view $E(f) -\area f(M) $ as a measure of nonconformality.
For the identity map, $E(\id_M) = \area(M)$ and
$\del_p E(\id_M) = \del_p\area(M)$,
confirming that $M$ is minimal if and only if $\id_M$ is harmonic. 

\section{Vector bundles on polyhedral manifolds}\index{vector bundle (discrete)}

We now give a general definition of vector bundles and
connections on polyhedral manifolds; this leads to another
interpretation of the \Gauss/ curvature for a polyhedral surface.

A polyhedral $n$-manifold~$P^n$ means a CW-complex
which is homeomorphic to an $n$-dimensional manifold,
and which is regular and satisfies the intersection condition.
(Compare~\cite{Ziegler-highg} in this volume.)  That is, each $n$-cell
(called a \emph{facet}) is embedded in~$P^n$ with no identifications on its
boundary, and the intersection of any two cells (of any dimension)
is a single cell (if nonempty).  Because~$P^n$ is topologically
a manifold, each $(n-1)$--cell
(called a \emph{ridge})\index{ridge (of complex)}
is contained in exactly two facets.\index{face (of polytope)}

\begin{defn}
A \emph{discrete vector bundle}~$V^k$ of rank~$k$ over~$P^n$
consists of a vectorspace $V_f\isom\E^k$ for each facet~$f$ of~$P$.
A \index{connection (discrete)}\emph{connection} on~$V^k$ is a choice of isomorphism~$\phi_r$
between~$V_f$ and~$V_{f'}$ for each ridge $r=f\cap f'$ of~$P$.
We are most interested in the case where the vectorspaces~$V_f$ have
inner products, and the isomorphisms~$\phi_r$ are orthogonal.
\end{defn}

Consider first the case $n=1$, where $P$ is a polygonal curve.
On an arc (an open curve) any vector bundle is trivial.
On a loop (a closed curve), a vector bundle of rank~$k$
is determined (up to isomorphism) simply by its \dfn{holonomy}
around the loop, an automorphism $\phi:\E^k\to\E^k$.

Now suppose $P^n$ is linearly immersed in~$\E^d$ for some~$d$.
That is, each $k$-face of~$P$ is mapped homeomorphically to
a convex polytope in an affine $k$-plane in~$\E^d$, and
the star of each vertex is embedded.
Then it is clear how to define the discrete tangent bundle~$T$
(of rank~$n$) and normal bundle~$N$ (of rank~$d-n$).
Namely, each $T_f$ is the $n$-plane parallel
to the affine hull of the facet~$f$,
and $N_f$ is the orthogonal $(d-n)$--plane.  These inherit
inner products from the euclidean structure of~$\E^d$.

There are also natural analogs of the \ix{Levi-Civita connection}s
on these bundles.  Namely, for each ridge~$r$, let $\alpha_r\in[0,\pi)$
be the exterior dihedral angle bewteen the facets~$f_i$ meeting along~$r$.
Then let $\phi_r:\E^d\to\E^d$ be the simple rotation by this angle,
fixing the affine hull of~$r$ (and the space orthogonal to the affine
hull of the~$f_i$).  We see that $\phi_r$ restricts to give maps
$T_f\to T_{f'}$ and $N_f\to N_{f'}$; these form the connections we want.
(Note that $T\oplus N=\E^d$ is a trivial vector bundle over~$P^n$,
but the maps $\phi_r$ give a nontrivial connection on it.)

Consider again the example of a closed polygonal curve $P^1\subset\E^d$.
The tangent bundle has rank~$1$ and trivial holonomy.  The holonomy of the
normal bundle is some rotation of~$\E^{d-1}$.  For $d=3$ this rotation
of the plane~$\E^2$ is specified by an angle equal (modulo~$2\pi$)
to the \ix{writhe} of~$P^1$.
(To define the writhe of a curve as a real number,
rather than just modulo~$2\pi$, requires a bit more care, and requires
the curve~$P$ to be embedded.)

Next consider a two-dimensional polyhedral surface~$P^2$ and its tangent bundle.
Around a vertex~$p$ we can compose the cycle of isomorphisms~$\phi_e$
across the edges incident to~$p$.  This gives a self-map $\phi_p:T_f\to T_f$.
This is a rotation of the tangent plane
by an angle which---it is easy to check---equals
the discrete \Gauss/ curvature\index{Gauss curvature!discrete}~$K_p$.

Now consider the general case of the tangent bundle to a polyhedral
manifold $P^n\subset\E^d$.  Suppose~$p$ is a codimension-two face of~$P$.
Then composing the ring of isomorphisms across the ridges incident to~$p$
gives an automorphism of~$\E^n$ which is the local holonomy, or \emph{curvature}
of the Levi-Civita connection around~$p$.  We see that this is a rotation
fixing the affine hull of~$p$.  To define this curvature (which can be
interpreted as a sectional curvature in the two-plane normal to~$p$)
as a real number and not just modulo~$2\pi$, we should look again at
the angle defect around~$p$, which is $2\pi-\sum\beta_i$ where the $\beta_i$
are the interior dihedral angles along~$p$ of the facets~$f$ incident to~$p$.

In the case of a hypersurface~$P^{d-1}$ in~$\E^d$, the one-dimensional
normal bundle is locally trivial: there is no curvature or local holonomy
around any~$p$.  Globally, the normal bundle is of course trivial exactly
when~$P$ is orientable.


\newcommand{\etalchar}[1]{$^{#1}$}
\providecommand{\bysame}{\leavevmode\hbox to3em{\hrulefill}\thinspace}

\end{document}